\documentstyle{amsppt}
\baselineskip18pt
\magnification=\magstep1
%\NoPageNumbers
%\NoRunningHeads
%\pagewidth{4.5in}
%\pageheight{7.0in}
\pagewidth{30pc}
\pageheight{45pc}
%\hoffset .5in
%\hsize 6in
%common header stuff

\hyphenation{co-deter-min-ant co-deter-min-ants pa-ra-met-rised
pre-print pro-pa-gat-ing pro-pa-gate
fel-low-ship Cox-et-er dis-trib-ut-ive}
\def\leaderfill{\leaders\hbox to 1em{\hss.\hss}\hfill}
\def\A{{\Cal A}}

\def\H{{\Cal H}}

\def\L{{\Cal L}}
\def\latl#1{{\Cal L}_L^{#1}}
\def\latr#1{{\Cal L}_R^{#1}}
\def\lat#1{{\Cal L}^{#1}}

\

\def\ldescent#1{{\Cal L (#1)}}
\def\rdescent#1{{\Cal R (#1)}}

\def\ti{\widetilde}

\def\d{{\delta}}

\def\e{{\varepsilon}}

\def\th{{\theta}}

\def\l{{\lambda}}

\def\T{{\widetilde T}}
\def\te{\widetilde t}

\def\bi{{\bold i}}

\def\bs{{\bold s}}
\def\bt{{\bold t}}
\def\bu{{\bold u}}

\def\b0{\text{\bf 0}}

\def\ra{{\ \longrightarrow \ }}

\def\coms{\text{\rm Co}(X, S)}

\def\lan{{\langle}}
\def\ran{{\rangle}}

\def\zed{{\Bbb Z}}

\def\enn{{\Bbb N}}

\def\boxit#1{\vbox{\hrule\hbox{\vrule \kern3pt
\vbox{\kern3pt\hbox{#1}\kern3pt}\kern3pt\vrule}\hrule}}
\def\rabbit{\vbox{\hbox{\kern0pt
\vbox{\kern0pt{\hbox{---}}\kern3.5pt}}}}

\def\tableau#1{
        \hbox {
                \hskip -10pt plus0pt minus0pt
                \raise\baselineskip\hbox{
                \offinterlineskip
                \hbox{#1}}
                \hskip0.25em
        }
}

\def\tabCol#1{
\hbox{\vtop{\hrule
\halign{\strut\vrule\hskip0.5em##\hskip0.5em\hfill\vrule\cr\lower0pt
\hbox\bgroup$#1$\egroup \cr}
\hrule
} } \hskip -10.5pt plus0pt minus0pt}

\def\CR{
        $\egroup\cr
        \noalign{\hrule}
        \lower0pt\hbox\bgroup$
}

% Set up the map arrows for commutative diagrams.
%\def\mapright#1{\smash{
%     \mathop{\longrightarrow}\limits^{#1}}}

%Set up macro for commutative diagrams etc. (see Ex. 18.46 in TeXbook)

\def\blank#1#2{
%\hbox to {{#1}}{\vbox to {{#2}}}
\hbox to #1{\hfill \vbox to #2{\vfill}}
}

%Ross's table macros

\def\strut{\vrule height10pt depth5pt width0pt}

%Paul's stuff

\def\ignore#1{{}}

\def\seca{1}
\def\secb{2}
\def\secc{3}
\def\secd{4}
\def\sece{5}
\def\secf{6}
\def\secconc{7}
\def\twola{2l+1}
\def\twolb{2l+2}
\def\twolc{2l+3}
\def\twold{2l+4}

\topmatter
\title Star reducible Coxeter groups
\endtitle

\author R.M. Green \endauthor
\affil Department of Mathematics \\ University of Colorado \\
Campus Box 395 \\ Boulder, CO  80309-0395 \\ USA \\ {\it  E-mail:}
rmg\@euclid.colorado.edu \\
\newline
\endaffil

\abstract 
We define ``star reducible'' Coxeter groups to be those Coxeter groups for
which every fully commutative element (in the sense of Stembridge) is
equivalent to a product of commuting generators by a sequence of 
length-decreasing star operations (in the sense of Lusztig).
We show that the Kazhdan--Lusztig bases of these groups have a nice
projection property to the Temperley--Lieb type quotient, and furthermore
that the images of the basis elements $C'_w$ (for fully commutative $w$)
in the quotient have structure constants in ${\Bbb Z}^{\geq 0}[v, v^{-1}]$.
We also classify the star reducible Coxeter groups and show that they form
nine infinite families with two exceptional cases.
\endabstract


\endtopmatter

\centerline{\bf To appear in the Glasgow Mathematical Journal}

\head Introduction \endhead

Let $(W, S)$ be a Coxeter group, with finite generating set $S$.  
Stembridge \cite{{\bf 19}} introduced the set $W_c$
of fully commutative elements of $W$ as those 
for which any two reduced expressions in the generators are equivalent 
via iterated application of short braid relations, that is, relations of 
the form $ss' = s's$, where $s, s' \in S$.  For example, if $w$ is a product
of commuting generators from $S$, then $w$ is fully commutative.

If $I = \{s, s'\} \subseteq S$ is a pair of noncommuting
Coxeter generators, then $I$ induces four partially defined maps from $W$
to itself, known as star operations.  A star operation, when it is defined,
respects the partition $W = W_c \dot\cup (W \backslash W_c)$ of the Coxeter
group, and increases or decreases the length of the element to which it is
applied by $1$.

In this paper we will analyse the situation where every fully commutative
element can be reduced to a product of commuting generators from $S$ by
iterated application of length-decreasing star operations; this property
is called ``Property F'' in \cite{{\bf 12}}, as it is essentially the same as
Fan's notion of cancellability in \cite{{\bf 5}}.  Groups with this property 
are the eponymous ``star reducible Coxeter groups'', and they include the
finite Coxeter groups as a subclass.

We shall show (Theorem \secd.1) that arbitrary elements of star reducible 
Coxeter 
groups have reduced expressions of a particularly nice type, which allows us 
to prove (Theorem \secd.3) a strong form of a certain conjectured projection 
property (in the sense of \cite{{\bf 14}, {\bf 18}}) for the associated 
Kazhdan--Lusztig basis
$\{C'_w: w \in W\}$.  This has some strong consequences (Theorem \secd.6) for
the Kazhdan--Lusztig type basis $\{c_w : w \in W_c\}$ introduced by J.
Losonczy and the author for a Temperley--Lieb type quotient of the Hecke
algebra $\H$ associated to $W$.  In the star reducible case, this basis
turns out simply to be the projection of the Kazhdan--Lusztig basis 
elements $\{C'_w : w \in W_c\}$.  Furthermore, there is a simple inductive
construction for the $c_w$, and the $c$-basis can be shown to have nonnegative
structure constants, that is, structure constants that are Laurent 
polynomials with nonnegative coefficients.  One of the reasons this is
interesting is that in many cases (see \cite{{\bf 12}, \S6} and \cite{{\bf 15},
Theorem 2.2.3, \S3.1}), these structure constants are also 
structure constants for the Kazhdan--Lusztig basis, whose positivity is
generally very difficult to prove.

Finally (Theorem \secf.3), we classify all star reducible Coxeter groups 
for which $S$ is a
finite set.  This class of groups contains the seven infinite families of 
groups ($A$, $B$, $D$, $E$, $F$, $H$ and $I$) for which $W_c$ is finite, 
which were classified independently by Graham \cite{{\bf 8}} and Stembridge 
\cite{{\bf 19}}, as well as three other infinite families (one of which 
subsumes type $I$) and two exceptional cases.  

Combining the main result of this paper (Theorem \secd.6) with the 
classification of star reducible Coxeter groups (Theorem \secf.3), one
obtains an extensive class of examples of situations
where the projection of the Kazhdan--Lusztig basis elements $C'_w$ (for
fully commutative $w$) to the Temperley--Lieb quotient have positive
structure constants.  These quotients are useful because they provide 
combinatorially tractable models for Kazhdan--Lusztig theory that are useful 
for formulating and checking conjectures, and in a future paper we plan to 
explain the application of the quotient algebras to the representation 
theory of the corresponding Lie algebras.
Our results here also provide unifying conceptual proofs for various results
already in the literature.

% 1. preliminaries

\head \seca. Preliminaries \endhead

Let $X$ be a Coxeter graph, of arbitrary type,
and let $W = W(X)$ be the associated Coxeter group with distinguished
(finite) set of generating involutions $S(X)$.  
(The reader is referred to
\cite{{\bf 1}} or \cite{{\bf 16}} for details of the theory of Coxeter groups.)
In other words, $W = W(X)$ is 
given by the presentation $$
W = \lan S(X) \ | \ (st)^{m(s, t)} = 1 \text{ for } m(s, t) < \infty \ran
,$$ where $m(s, s) = 1$ and $m(s, t) = m(t, s)$.
It turns out that the elements 
of $S = S(X)$ are distinct as group elements, and that $m(s, t)$
is the order of $st$.

Denote by $S^*$ the
free monoid on $S = S(X)$. We call the elements of $S$ {\it letters} and
those of $S^*$ {\it words}. The {\it length} of a word is the
number of factors required to write the word as a product of
letters. Let $\phi : S^* \longrightarrow W$ be the surjective
morphism of monoid structures satisfying $\phi(i) = s_i$ for all
$i\in S$. A word $\bi \in S^*$ is said to {\it represent} its
image $w=\phi(\bi)\in W$; furthermore, if the length of $\bi$ is
minimal among the lengths of all the words that represent $w$,
then we call $\bi$ a {\it reduced expression} for $w$. The {\it
length} of $w$, denoted by $\ell(w)$, is then equal to the length
of $\bi$.
A product $w_1w_2\cdots w_n$ of elements $w_i\in W$ is called
{\it reduced} if $\ell(w_1w_2\cdots w_n)=\sum_i\ell(w_i)$.  
We write $$
\ldescent{w} = \{s \in S : \ell(sw) < \ell(w)\}
$$ and $$
\rdescent{w} = \{s \in S : \ell(ws) < \ell(w)\}
.$$  The set $\ldescent{w}$ (respectively, $\rdescent{w}$) is called the 
{\it left} (respectively, {\it right}) {\it descent set} of $w$.

The {\it commutation monoid} $\coms$ is the quotient of the free
monoid $S^*$ by the congruence $\equiv$ generated by the commutation
relations: $$
st \equiv ts \text{ for all } s, t \in S \text{ with } \phi(s)\phi(t) = 
\phi(t)\phi(s)
;$$ note that, as a monoid, $W$ is a quotient of $\coms$.

The elements of $\coms$, which computer scientists call {\it traces}
\cite{{\bf 3}}, have the following normal form, often called the Cartier--Foata
normal form (see \cite{{\bf 2}}).

\proclaim{Theorem \seca.1 (Cartier--Foata normal form)}
Let $\bs$ be an element of the commutation monoid $\coms$.  Then $\bs$ has
a unique factorization in $\coms$ of the form $$
\bs = \bs_1 \bs_2 \cdots \bs_p
$$ such that each $\bs_i$ is a product of distinct commuting elements of
$S$, and such that for each $1 \leq j < p$ and each generator $t \in S$
occurring in $\bs_{j+1}$, there is a generator $s \in S$ occurring in
$\bs_j$ such that $st \ne ts$.
\qed\endproclaim

\remark{Remark \seca.2}
The Cartier--Foata normal form may be defined inductively, as follows.  
If we define 
$\ldescent{\bs}$ to be the set of possible first letters in all the words
$\bs'$ for which $\bs' \equiv \bs$ in $\coms$, then $\bs_1$ is just the 
product of
the elements in $\ldescent{\bs}$.  Since $\coms$ is a cancellative monoid,
there is a unique element $\bs' \in \coms $ with $\bs = \bs_1 \bs'$.  If $$
\bs' = \bs_2 \cdots \bs_p
$$ is the Cartier--Foata normal form of $\bs'$, then $$
\bs_1 \bs_2 \cdots \bs_p
$$ is the Cartier--Foata normal form of $\bs$.
\endremark

Denote by $\H = \H(X)$ the Hecke
algebra associated to $W$.  This is a $\zed[q, q^{-1}]$-algebra 
with a basis consisting of (invertible) elements $T_w$, with $w$ ranging over 
$W$, satisfying $$T_s T_w = 
\cases
T_{sw} & \text{ if } \ell(sw) > \ell(w),\cr
q T_{sw} + (q-1) T_w & \text{ if } \ell(sw) < \ell(w),\cr
\endcases$$ where $\ell$ is the length function on the Coxeter group
$W$, $w \in W$, and $s \in S$.

For many applications it is convenient to introduce an $\A$-form of $\H$,
where $\A = \zed[v, v^{-1}]$ and $v^2 = q$, and to define a scaled version
of the $T$-basis, $\{\T_w : w \in W\}$, where $\T_w := v^{-\ell(w)} T_w$.
Unless otherwise stated, we will use the $\A$-form of $\H$ from now on, and
we will denote the $\zed[q, q^{-1}]$-form by $\H_q$.
We will write $\A^+$ and $\A^-$ for $\zed[v]$ and $\zed[v^{-1}]$, respectively,
and we denote the $\zed$-linear ring
homomorphism $\A \ra \A$ exchanging $v$ and $v^{-1}$ by $\bar{\ }$.  We
can extend $\bar{\ }$ to a ring automorphism of $\H$
(as in \cite{{\bf 7}, Theorem 11.1.10}) by the condition that $$
\overline{\sum_{w \in W} a_w \T_w} := \sum_{w \in W} \overline{a_w}
\T_{w^{-1}}^{-1}
,$$ where the $a_w$ are elements of $\A$.

In \cite{{\bf 17}}, Kazhdan and Lusztig proved the following

\proclaim{Theorem \seca.3. (Kazhdan, Lusztig)}
For each $w \in W$, there exists a unique $C'_w \in \H$ such that both
$\overline{C'_w} = C'_w$ and $$
C'_w = \T_w + \sum_{y < w} a_y \T_y
,$$ where $<$ is the Bruhat order on $W$ and $a_y \in v^{-1} \A^-$.
The set $\{C'_w : w \in W\}$ forms an $\A$-basis for $\H$.
\qed\endproclaim

Following \cite{{\bf 7}, \S11.1}, we denote the coefficient of $\T_y$ in $C'_w$ by 
$P^*_{y, w}$.  The {\it Kazhdan--Lusztig polynomial} $P_{y, w}$ is then 
given by $v^{\ell(w) - \ell(y)} P^*_{y, w}$.

Let $J(X)$ be the two-sided ideal of $\H$ generated by the elements $$
\sum_{w \in \lan s, s' \ran} T_w,
$$ where $(s, s')$ runs over all pairs of elements of $S$
that correspond to adjacent nodes in the Coxeter graph, and $\lan s, s' \ran$
is the parabolic subgroup generated by $s$ and $s'$.   
(If the nodes corresponding to $(s, s')$ are connected by a
bond of infinite strength, then we omit the corresponding relation.)

Following Graham \cite{{\bf 8}, Definition 6.1}, we define the {\it generalized
Temperley--Lieb algebra} $TL(X)$ to be
the quotient $\A$-algebra $\H(X)/J(X)$.  We denote the corresponding
epimorphism of algebras by $\th : \H(X) \ra TL(X)$.  Since the generators
of $J(X)$ lie in $\H_q(X)$, we also obtain a $\zed[q, q^{-1}]$-form
$TL_q(X)$, of $TL(X)$.
Let $t_w$ (respectively, $\te_w$) denote the image in 
$TL(X)$ of the basis element $T_w$ (respectively, $\T_w$) of $\H$.

Call an element $w \in W$ {\it complex} if it can be written 
as a reduced product $x_1 w_{ss'} x_2$, where $x_1, x_2 \in W$ and
$w_{ss'}$ is the longest element of some rank 2 parabolic subgroup 
$\lan s, s'\ran$ such that $s$ and $s'$ correspond to adjacent nodes
in the Coxeter graph.
An element $w \in W$ is said to be {\it weakly complex} if it is
complex and of the form $w = su$, where $u$ is not complex and $s \in S$.
In this case, we must have $su > u$.

Denote by $W_c(X)$ the set of all elements of $W$
that are not complex.  The elements of $W_c$ are the {\it fully commutative}
elements of \cite{{\bf 19}}; they are characterized by the property that any two
of their reduced expressions may be obtained from each other by repeated
commutation of adjacent generators; in other words, all reduced expressions
are equal as elements of $\coms$.  Each reduced expression for $w$ has a 
Cartier--Foata normal form, by considering it as an element of $\coms$, and
this normal form is an invariant of $w$ if and only if $w$ is fully
commutative.

We define the $\A^-$-submodule $\L$ of
$TL(X)$ to be that generated by $\{\te_w : w \in W_c\}$.  We define
$\pi : \L \ra \L/v^{-1}\L$ to be the canonical $\zed$-linear projection.

By \cite{{\bf 13}, Lemma 1.4}, the ideal $J(X)$ is fixed by $\bar{\ }$, so 
$\bar{\ }$ induces an involution on $TL(X)$, which we also denote by 
$\bar{\ }$.

The following result is an analogue of Theorem \seca.3 for the quotient 
algebra.

\proclaim{Theorem \seca.4}
\item{\rm (i)}
{The set $\{t_w : w \in W_c\}$ is a $\zed[q, q^{-1}]$-basis for $TL_q(X)$.
The set $\{\te_w : w \in W_c\}$ is an $\A$-basis for $TL(X)$, and an
$\A^-$-basis for $\L$.}
\item{\rm (ii)}
{For each $w \in W_c$, there exists a unique $c_w \in TL(X)$ such that both
$\overline{c_w} = c_w$ and $\pi(c_w) = \pi(\te_w)$.  Furthermore, we have $$
c_w = \te_w + \sum_{{y < w} \atop {y \in W_c}} a_y \te_y
,$$ where $<$ is the Bruhat order on $W$, and $a_y \in \A^-$ for all $y$.}
\item{\rm (iii)}
{The set $\{c_w : w \in W_c\}$ forms an $\A$-basis for $TL(X)$ and an
$\A^-$-basis for $\L$.}
\item{\rm (iv)}
{If $x \in \L$ and $\bar{x} = x$, then $x$ is a
$\zed$-linear combination of the $c_w$.}
%\item{\rm (v)}
%{There is an $\A$-linear anti-automorphism, $*$, of $TL(X)$ that sends
%$\te_w$ to $\te_{w^{-1}}$ and $c_w$ to $c_{w^{-1}}$ for all $w \in W_c$.}
\endproclaim

\demo{Proof}
This is a subset of \cite{{\bf 12}, Theorem 2.1}.  (Note that (i) is
due to Graham \cite{{\bf 8}, Theorem 6.2}, and (ii) and (iii) are essentially
due to J. Losonczy and the author \cite{{\bf 13}, Theorem 2.3}.)
\qed\enddemo

Let $W$ be any Coxeter group and let $I = \{s, t\} \subseteq S$ be a
pair of noncommuting generators whose product has order $m$ (where $m = \infty$
is allowed).  Let $W^I$ denote the set of all $w \in W$ satisfying 
$\ldescent{w} \cap I = \emptyset$.  Standard properties of Coxeter groups 
\cite{{\bf 16}, \S5.12} show that any element $w \in W$ may be uniquely
written as $w = w_I w^I$, where $w_I \in W_I = \lan s, t \ran$
and $\ell(w) = \ell(w_I) + \ell(w^I)$.  
There are four possibilities for elements $w \in W$:
\item{(i)}{$w$ is the shortest element in the coset $W_I w$, so $w_I = 1$ and 
$w \in W^I$;}
\item{(ii)}{$w$ is the longest element in the coset $W_I w$, so $w_I$ is 
the longest element of $W_I$ (which can only happen if $W_I$ is finite);}
\item{(iii)}{$w$ is one of the $(m-1)$ elements $sw^I$, $tsw^I$, $stsw^I, 
\ldots$;}
\item{(iv)}{$w$ is one of the $(m-1)$ elements $tw^I$, $stw^I$, $tstw^I, 
\ldots$.}

The sequences appearing in (iii) and (iv) are called {\it (left) 
$\{s, t\}$-strings},
or {\it strings} if the context is clear.  If $x$ and $y$ are two elements 
of an $\{s, t\}$-string such that $\ell(x) = \ell(y) - 1$, we call the pair 
$\{x, y\}$ {\it left $\{s, t\}$-adjacent}, and we say that $y$ is 
{\it left star reducible to $x$}.

The above concepts all have right-handed counterparts, leading to the notion of
{\it right $\{s, t\}$-adjacent} and {\it right star reducible} pairs of 
elements, and coset decompositions $({^Iw})({_Iw})$.

If there is a (possibly trivial) sequence $$
x = w_0, w_1, \ldots, w_k = y
$$ where, for each $0 \leq i < k$, $w_{i+1}$ is left star reducible or
right star reducible to $w_i$ 
with respect to some pair $\{s_i, t_i\}$, we 
say that $y$ is {\it star reducible to $x$}.  Because star reducibility
decreases length, it is clear that this defines a partial order on $W$.

If $w$ is an element of an $\{s, t\}$-string, $S_w$, we have 
$\{\ell(sw), \ell(tw)\}$ = $\{\ell(w) - 1, 
\ell(w) + 1\}$; let us assume without loss of generality that $sw$ is longer
than $w$ and $tw$ is shorter.  If $sw$ is an element of $S_w$, we define
$^*w = sw$; if not, $^*w$ is undefined.  If $tw$ is an element of $S_w$,
we define $_*w = tw$; if not, $_*w$ is undefined.

There are also obvious right handed analogues to the above concepts,
so the symbols $w^*$ and $w_*$ may be used with the analogous meanings.

\example{Example \seca.5}
In the Coxeter group of type $B_2$ with $w = ts$, we have $$
{_*w} = s, \ {^*w} = sts, \ w_* = t \text{\ and\ } w^* = tst
.$$  If $x = sts$ then $^*x$ and $x^*$ are undefined; if $x = t$ then
$_*x$ and $x_*$ are undefined.
\endexample

\definition{Definition \seca.6}
We say that a Coxeter group $W(X)$, or its Coxeter graph $X$, is {\it star
reducible} if every element of $W_c$ is star reducible to a
product of commuting generators from $S$.
\enddefinition

% 2. weak Property W (uses heaps)

\head \secb. Acyclic monomials \endhead

In order to derive some of the results in this paper, and \S\secb\  in 
particular, we will need to use the author's theory of acyclic heaps 
\cite{{\bf 10}, {\bf 11}}.
Heaps, as introduced by Viennot in \cite{{\bf 21}}, are certain combinatorial
structures associated to elements of $\coms$; they are known as ``dependence
graphs'' in the computer science literature \cite{{\bf 3}}.
However, in order to keep the paper as accessible as possible, we will avoid 
mention of heaps and work directly with monomials, or traces.  All Coxeter
groups in \S\secb\  will be star reducible.

\proclaim{Theorem \secb.1}
Let $(W, S)$ be a star reducible Coxeter group.  There is a unique 
function $h : \coms \ra \zed^{\geq 0}$ with the 
following properties.
\item{\rm (i)}
{If $\bu \in \coms$ and $s, t \in S$ are noncommuting generators, then $
h(st\bu) = h(t\bu)
$ and $
h(\bu ts) = h(\bu t)
.$}
\item{\rm (ii)}
{If $\bu \in \coms$ is represented by a monomial $$
s_1 s_2 \cdots s_r
$$ that is a reduced expression for some $w \in W_c$, then $h(\bu) = 0$.}
\item{\rm (iii)}
{If $\bu = \bu_1 s s \bu_2$ for some generator $s \in S$,
and $\bu' = \bu_1 s \bu_2$, then $h(\bu) = h(\bu') + 1$.}
\item{\rm (iv)}
{If $\bu = \bu_1 s t s \bu_2$ for some noncommuting generators $s, t \in S$,
and $\bu' = \bu_1 s \bu_2$, then $h(\bu) = h(\bu')$.}
\item{\rm (v)}
{If $\bu = \bu_1 s \bu_2$ for some generator $s \in S$,
and $\bu' = \bu_1 \bu_2$, then $|h(\bu) - h(\bu')| \leq 1$.}
\endproclaim

\demo{Proof}
Let $k$ be a field.

According to \cite{{\bf 21}, Proposition 3.4}, elements $\bu$ of $\coms$ are in 
bijection with certain heaps $[E, \leq, \e]$ (see \cite{{\bf 10}}, and 
\cite{{\bf 10}, Proposition 3.1.4} in particular, for more details on these 
concepts and the notation).
Let $h(\bu) = \dim H_1(E, k)$; it will turn out that the definition is
independent of $k$.  

Part (i) follows from the proof of the inductive
step in \cite{{\bf 10}, Proposition 2.2.3}.

Since $W$ is star reducible, it follows from using (i) repeatedly that 
(ii) is true if and only if it is true when $\bu$ is a product of distinct 
commuting generators.  In this case, the claim follows from the proof of
the base case of the induction in \cite{{\bf 10}, Proposition 2.2.3}.

Part (iii) is a restatement of \cite{{\bf 10}, Lemma 2.3.4}, part (iv) is
a restatement of \cite{{\bf 10}, Lemma 2.3.5}, and part (v) is a restatement of
\cite{{\bf 10}, Theorem 2.1.1} (star reducibility plays no role
in these proofs).

It follows from \cite{{\bf 1}, Theorem 3.3.1 (i)} that
the elements of $\coms$ corresponding to reduced expressions of some
$w \in W$ are precisely those that have no monomial representative of
the form $\bu_1 s s \bu_3$, where $s \in S$, and no monomial representative
$\bu_1 \bu_2 \bu_3$ where $\bu_2$ is an alternating product of
$m(s, t) > 2$ occurrences of $s$ and $t$.  It follows from this that any
element of $\coms$ can be transformed into an element of $\coms$ corresponding
to a reduced expression for some $w \in W_c$ by repeatedly applying 
transformations of the form $ss \mapsto s$ or $sts \mapsto s$, as used
in parts (iii) and (iv).  Applying (ii), we see there is at most one
function $h$ satisfying (ii), (iii) and (iv).  This proves uniqueness of
$h$ and also shows that the definition is independent of the choice of field
$k$.
\qed\enddemo

\definition{Definition \secb.2}
In the set-up of Theorem \secb.1, an element $\bu$ of $\coms$ (and, by 
extension, an element of $S^*$ representing $\bu$) is called an
{\it acyclic monomial} if $h(\bu) = 0$.  (The acyclic monomials are those
that correspond to the acyclic heaps of \cite{{\bf 10}, {\bf 11}}.)
\enddefinition

For our purposes in this paper, it is convenient to work with another
basis of $TL(X)$, namely the monomial basis.  Although the fact that this is
a basis is well-known, we provide a proof since there does not seem to be an 
easily available general proof in the literature.

\definition{Definition \secb.3}
Let $W$ be a Coxeter group and let $w \in W_c$ be a fully commutative
element.  Let $$
w = s_1 s_2 \cdots s_r
$$ be a reduced expression for $w$.  For each $s \in S$, let
$b_s = v^{-1} \te_1 + \te_s$, then define $b_w \in TL(X)$ by $$
b_w := b_{s_1} b_{s_2} \cdots b_{s_r}
.$$
\enddefinition

Note that the element $b_w$ is well-defined precisely because any
two reduced expressions for $w$ are commutation equivalent.

\proclaim{Proposition \secb.4}
The set $\{ b_w : w \in W_c\}$ is a free $\A$-basis for $TL(X)$, and
$\overline{b_w} = b_w$ for all $w \in W_c$.
\endproclaim

\demo{Proof}
The second assertion follows from the fact that $\bar{\ }$ is a ring 
endomorphism of $TL(X)$ that fixes the generators $b_s = c_s (s \in S)$.

To prove the first assertion, first observe that 
by definition of the ideal $J(X)$, we have the relation $$\eqalignno{
\te_{w_{ss'}} &= - \sum_{w  \in \lan s, s' \ran, w < w_{ss'}}
v^{\ell(w) - \ell(w_{ss'})} \te_w. & (1)
}$$ in $TL(X)$, where $w_{ss'}$ is the longest element in the parabolic
subgroup $\lan s, s' \ran$ of $W$.  This has the consequence that 
any monomial $$
\te_{s_1} \te_{s_2} \cdots \te_{s_k}
,$$ where all $s_i \in S$, can be expressed as a linear combination of
basis elements $\te_x$ for which $\ell(x) \leq k$.  Now let $x \in W_c$ and
let $s_1 s_2 \cdots s_r$ be a reduced expression for $x$.  Since $$
b_x = b_{s_1} b_{s_2} \cdots b_{s_r}
,$$ we have $$
b_x = 
(v^{-1} \te_1 + \te_{s_1})
(v^{-1} \te_1 + \te_{s_2})
\cdots
(v^{-1} \te_1 + \te_{s_r})
.$$  Expanding the parentheses and using equation (1), we see that $$
b_x = \te_x + \sum_{{y \in W_c} \atop {\ell(y) < \ell(x)}} a_y \te_y
$$ for some coefficients $a_y \in \A$.  It is now clear that the set in
the statement is a basis, and that the change of basis matrix from the
$\te$-basis to the $b$-basis is unitriangular.
\qed\enddemo

It will be convenient to have a presentation of $TL(X)$ in terms of the
generators $b_s$; compare with \cite{{\bf 8}, Proposition 9.5}.

\definition{Definition \secb.5}
We define the Chebyshev polynomials of the second kind
to be the elements of $\zed[x]$ given by the conditions $P_0(x) = 1$,
$P_1(x) = x$ and $$
P_n(x) = xP_{n-1}(x) - P_{n-2}(x)
$$ for $n \geq 2$.  
If $f(x) \in \zed[x]$, 
we define $f_b^{s, t}(x)$ to be 
the element of $TL(X)$ given by the linear extension of the map sending
$x^n$ to the product $$
\underbrace{b_s b_t \ldots}_{n \text{ factors}}
$$ of alternating factors starting with $b_s$.  
\enddefinition

\proclaim{Proposition \secb.6}
As a unital $\A$-algebra, $TL(X)$ is given by generators $\{b_s : s \in S\}$
and relations $$\eqalignno{
b_s^2 &= \d b_s, & (2)\cr
b_s b_t &= b_t b_s \text{\quad if } m(s, t) = 2, & (3)\cr
(x P_{m-1})_b^{s, t}(x) &= 0 \text{\quad if } 2 < m = m(s, t) < \infty, 
& (4)\cr
}$$ where $\d := (v + v^{-1})$.
\endproclaim

\demo{Proof}
This follows from \cite{{\bf 12}, Corollary 6.5} and its proof, which shows that
if $2 < m(s, t) < \infty$,  $$
(x P_{m-1})_b^{s, t}(x)
$$ is the image in $TL(X)$ of $C'_{w_{st}}$.  (A similar result appears in 
\cite{{\bf 8}, Proposition 9.5}.)
\qed\enddemo

\example{Example \secb.7}
Relation (4) reads $$\eqalign{
b_s b_t b_s - b_s &= 0 \text{\quad if } m = 3,\cr
b_s b_t b_s b_t - 2 b_s b_t &= 0 \text{\quad if } m = 4,\cr
b_s b_t b_s b_t b_s - 3 b_s b_t b_s + b_s &= 0 \text{\quad if } m = 5, 
\text{ and}\cr
b_s b_t b_s b_t b_s b_t - 4 b_s b_t b_s b_t + 3 b_s b_t &= 0
\text{\quad if } m = 6.\cr
}$$
\endexample

\remark{Remark \secb.8}
Since the relations (3) all occur in $\coms$, it makes sense, given an
element $\bs \in \coms$ represented by a monomial $s_1 s_2 \cdots s_r$, 
to define an element $b(\bs) \in TL(X)$ by $$
b(\bs) := b_{s_1} b_{s_2} \cdots b_{s_r}
.$$
\endremark

The following lemma is the generalization of \cite{{\bf 10}, Theorem 3.2.3}
alluded to in \cite{{\bf 10}, \S4.1}.

\proclaim{Lemma \secb.9}
Let $(W, S)$ be a star reducible Coxeter group, let $s_1 s_2 \cdots s_r$ be
an arbitrary monomial in $S^*$ representing the trace $\bs \in \coms$,
and let $b(\bs)$ be the element of $TL(X)$ given in Remark \secb.8.
Express $b$ as a linear combination of the monomial basis, namely $$
b(\bs) = \sum_{w \in W_c} \l_w b_w
.$$  Then $\l_w$ is an integer multiple of $\d^{h(\bs)}$, where $h$ is
as in Theorem \secb.1 and $\d = (v + v^{-1})$.
\endproclaim

\demo{Proof}
We claim that $TL(X)$ has the structure of a graded $\zed$-module $$
\bigoplus_{k \geq 0} M_k
,$$ where $M_k$ is the free $\zed$-module on the set $$
\{ \d^p b(\bt) \text{ such that } p \geq 0, \ \bt \in \coms 
\text{ and } p + h(\bt) = k \}
.$$  

The only nontrivial thing to check is that the grading is respected by
the relations of Proposition 2.6.  Relation (3) clearly respects the grading,
because it is a relation in $\coms$.  
Relation (2) respects the grading by Theorem \secb.1 (iii).

Note that relation (4) is a linear combination of monomials, each of which
can be transformed into any of the others by iterated substitutions of the
form $b_s b_t b_s \leftrightarrow b_s$ (see Example \secb.7 for clarification).
Although these substitutions are not generally valid relations in $TL(X)$,
it now follows from Theorem \secb.1 (iv) that relation (4) respects the
grading given.

Now consider the monomial $b(\bs)$.  By applying relations (2), (3) and (4)
repeatedly to express $b$ in terms of shorter monomials, we can write $b$ 
as a linear combination $$
b(\bs) = \sum_{w \in W_c} \l_w b_w
,$$ where $\l_w = n_w \d^{d_w}$ for some integer $n_w$ and nonnegative
integer $d_w$.  By Theorem \secb.1 (ii), all the monomials $b_w$ in the
sum are of the form $b(\bu)$, where $h(\bu) = 0$.  Since each
side of the equation lies in $M_{h(\bs)}$, it follows that $d_w = h(\bs)$,
as required.
\qed\enddemo

\proclaim{Lemma \secb.10}
If $W$ is a star reducible Coxeter group, then the $b$-basis and the
$c$-basis of $TL(X)$ have the same $\zed$-span.  In particular, the
$b$-basis is an $\A^-$-basis for $\L$.
\endproclaim

\demo{Proof}
Let $\bs = s_1 s_2 \cdots s_r$ be a reduced expression for $w \in W_c$, 
and write $$\eqalign{
\te_w &= \te_{s_1} \te_{s_2} \cdots \te_{s_r} \cr
&= 
(b_{s_1} - v^{-1})
(b_{s_2} - v^{-1})
\cdots
(b_{s_r} - v^{-1}).
}$$  Expanding the parentheses, we express $\te_w$ as a linear combination of
elements $(-v)^{-k} b(\bu)$, where $\bu$ is obtained from $\bs$ by deletion of
$k$ generators.  By Theorem \secb.1 (ii), $\bs$ is acyclic, so by Theorem
\secb.1 (v), we must have $h(\bu) \leq k$.  By Lemma \secb.9, if we express
$(-v)^{-k} b(\bu)$ in terms of the monomial basis, namely $$
(-v)^{-k} b(\bu) = \sum_{w \in W_c} (-v)^{-k} \l_w b_w
,$$ we see that $(-v)^{-k} \l_w \in \A^-$.

It follows from this that $\te_w$ is an $\A^-$-linear combination of monomial
basis elements.  Since any monomial in the $b_s$ is a linear combination of
basis monomials of shorter length, the above argument shows that the 
coefficient of $b_w$ in $\te_w$ is $1$.  This means that the change of basis
matrix from the $\te$-basis to the $b$-basis is unitriangular with entries
in $\A^-$ with respect to a suitable total ordering, and hence the inverse
of this matrix has the same properties, in other words, the monomial basis
elements lie in $\L$.

By Proposition \secb.4, $\overline{b_w} = b_w$ for any $w \in W_c$.  By
Theorem \seca.4 (iv), $b_w$ is a $\zed$-linear combination of $c$-basis 
elements.  By the above paragraph, we have $$
b_w = \te_w + \sum_{{x \in W_c} \atop {x < w}} \nu_x \te_x
$$ for certain $\nu_x \in \A^-$.  Applying $\pi$ to both sides and appealing
to Theorem \seca.4 (ii) and (iv), we have $$
b_w = c_w + \sum_{{x \in W_c} \atop {x < w}} \xi_x c_x
$$ for certain integers $\xi_x$.  This shows that the change of basis matrix
between the $b$-basis and the $c$-basis is unitriangular with entries in 
$\zed$ with respect to a suitable total ordering, from which it follows that
the $b$-basis and the $c$-basis have the same $\zed$-span.  This implies
that they also have the same $\A^-$-span, namely $\L$.
\qed\enddemo

% 3. Monomials and weakly complex elements

\head \secc.  Monomials and weakly complex elements \endhead

In \S\secc, we develop the properties of the lattice $\L$ by using the
monomial basis which, as we know from Lemma \secb.10, is an $\A^-$-basis 
for $\L$.

\proclaim{Lemma \secc.1}
Let $W$ be a star reducible Coxeter group.  Then, for $s \in S$, 
the set $$
\{ x \in TL(X) : b_s x = (v + v^{-1}) x \}
$$ is the free $\A$-submodule of $TL(X)$ with basis 
$B_s := \{b_y : y \in W_c, \ sy < y\}$.
\endproclaim

\demo{Proof}
If $y \in W_c$ is such that $sy < y$, it is clear that
$b_s b_y = \d b_s$ by relation (2), and it follows that the set $B_s$ is
contained in the required subset of $TL(X)$.

To finish the proof, it is enough to show that if $b(\bu) \in TL(X)$, then 
$b_s b(\bu)$ is a linear combination of elements $b_y$ with $y \in W_c$ and
$sy < y$.  

Let us say that a monomial $\bs = s_1 s_2 \cdots s_r \in S^*$ is 
``$s$-minimal'' if the following conditions are satisfied:
\item{1.}
{$s_i = s$ for some $1 \leq i \leq r$;}
\item{2.}
{$s_h \ne s$ for any $1 \leq h < i$;}
\item{3.}
{$s_h$ and $s_i$ commute for any $1 \leq h < i$.}

Condition 3 above means that it also makes sense to speak of an element 
$\bs \in \coms$ being $s$-minimal.

We see that applying one of the relations (2), (3) or
(4) to $b(\bs)$ results in a linear combination of monomials $b(\bt)$ where
$\bt$ is also $s$-minimal.  Repeating this argument shows that if
$b(\bs)$ is $s$-minimal, then it is a linear combination of $s$-minimal basis
elements.  However, the $s$-minimal basis elements are precisely those
basis elements $b_y$ where $y$ has a reduced expression beginning with $s$,
which implies by relation (2) that $b_s b_y = \d b_y$.  Since any monomial
of the form $b_s b(\bu)$ is $s$-minimal, the proof is complete.
\qed\enddemo

\remark{Remark \secc.2}
It is tempting to think from Lemma \secc.1 that if $b_y$ is a monomial basis
element such that $b_s b_y$ is $\d$ times another basis element, then $sy < y$,
but this is not true.  If $W$ is the (star reducible) Coxeter group
of type $B_3$, and $S = \{s_1, s_2, s_3\}$ is indexed so that
$m(s_1, s_2) = 4$ and $m(s_2, s_3) = 3$, then setting $y = s_1 s_2 s_1 s_3
\in W_c$ we have $$
b_{s_3} b_y = \d b_z
,$$ where $z = s_1 s_3 \in W_c$, even though $s_3 y > y$.  The $c$-basis
does not have this disadvantage, as will be clear from Theorem \secd.6 (ii)
below.
\endremark

We recall the following definition from \cite{{\bf 12}, \S4}.

\definition{Definition \secc.3}
Let $W' \subset W_c$.  We define $\lat{W'}$ to be the free $\A^-$-module with
basis $$
\{ \te_w : w \in W' \} \cup \{ v^{-1} \te_w : w \in W_c \backslash W' \}
.$$  
If $s, t \in S$ are noncommuting generators, 
$W_1 = \{ w \in W_c : sw < w \}$ and $W_2 = \{ w \in W_c : w = stu
\text{ reduced} \}$, we write
$\latl{s}$ and $\latl{st}$ for $\lat{W_1}$ and $\lat{W_2}$, respectively.

One can also define right handed versions, $\latr{s}$ and
$\latr{ts}$, of the above concepts, and of Lemma \secc.1.
\enddefinition

\proclaim{Lemma \secc.4}
Let $W$ be a star reducible Coxeter group.  Then the set $$
\{ b_y : y \in W_c, \ sy < y \} \cup \{ v^{-1} b_z : z \in W_c, \ sz > z \}
$$ is an $\A^-$-basis for $\latl{s}$.
\endproclaim

\demo{Proof}
Since the monomial basis is an $\A^-$-basis for $\L$ and there is a natural
bijection between the set in the statement and the defining $\A^-$-basis for
$\latl{s}$, the claim will follow if we can show that whenever we have 
$y \in W_c$ with $sy < y$, then $$
\pi(\te_y) = \pi(b_y) + \sum_{{w < y} \atop {sw < w}} 
\xi_w \pi(b_w)
,$$  where $w \in W_c$ in the sum and $\xi_w \in \zed$. 
Apart from the assertion that $sw < w$, this follows from the observations
relating the $b$-basis to the $\te$-basis made in the proof of Lemma \secb.10.

Since $y = sy'$ is reduced, we have $$
\te_y = \te_s \te_{y'} = (b_s - v^{-1}) \te_{y'}
,$$ and clearly $v^{-1} \te_{y'} \in v^{-1} \L$.  Since $\te_y \in \L$, it
follows that $b_s \te_{y'} \in \L$.  However, by Lemma \secc.1, we have $$
b_s \te_{y'} = \sum_{{w \leq y} \atop {sw < w}} \l_w b_w
,$$ where the sum is over $w \in W_c$ and we have $\l_w \in \A^-$ by 
Lemma \secb.10.  Since $\pi(\te_y) = \pi(b_s \te_{y'})$, the assertion follows.
\qed\enddemo

\proclaim{Lemma \secc.5}
Let $W$ be a star reducible Coxeter group, and let $s, t \in S$ be
noncommuting generators.  Then $
b_s \latl{t} \subseteq \latl{s}
$ and $
\te_s \latl{t} \subseteq \latl{s}
.$
\endproclaim

\demo{Proof}
The second assertion is immediate from the first and the identity
$b_s = (v^{-1} \te_1 + \te_s)$, so we concentrate on the first assertion.

Suppose that $y \in W_c$ is such that $ty < y$, and write $b_y = b(\bu)$
in the usual way, where $\bu \in \coms$.  By Theorem \secb.1 (ii), 
$h(\bu) = 0$, and by Theorem \secb.1 (i), $h(s \bu) = 0$ too.  Lemma \secb.9
now shows that $b_s b_y$ is a $\zed$-linear combination of basis elements
$b_w$, and then lemmas \secb.10 and \secc.1 show that $b_s b_y \in \latl{s}$.

Suppose now that $y \in W_c$ is such that $ty > y$, and write $b_y = b(\bu)$
as before.  In this case, $h(\bu) = 0$, and Theorem \secb.1 (v) shows that
$h(s \bu) \leq 1$.  Lemma \secb.9 then shows that $b_s b_y$ is a 
$v \A^-$-linear combination of basis elements $b_w$.  Lemmas \secb.10 and
\secc.1 show that $b_s b_y \in v \latl{s}$.

An application of Lemma \secc.4, combining the above two observations,
completes the proof.
\qed\enddemo

To prove the main result of \S\secc, we need to recall some of the 
combinatorial properties of weakly complex elements from \cite{{\bf 12}}.  The next
result shows that weakly complex elements respect the left and right weak
Bruhat orders.

\proclaim{Lemma \secc.6}
Let $W$ be any Coxeter group and let $w \in W_c$ be such that $sw
\not\in W_c$ for some $s \in S$.
If $u \in S$ and $y \in W$ are such that we have either $w = uy$ or $w = yu$
reduced, then either $sy \in W_c$ or $sy$ is weakly complex.
\endproclaim

\demo{Proof}
See \cite{{\bf 12}, Lemma 4.5 (iii)}.
\qed\enddemo

\proclaim{Lemma \secc.7}
Let $W$ be a star reducible Coxeter group, let $w \in W_c$ and 
$x = sw > w$, where $s \in S$.
Then one of the following situations must occur:
\item{\rm (i)}{$x$ is a product of commuting generators;}
\item{\rm (ii)}{$x \in W_c$ and there exists $I = \{s, t\} \subseteq S$ 
with $st \ne ts$ such that when $x = x_I x^I$, we have $\ell(x_I) > 1$;}
\item{\rm (iii)}{$x$ is weakly complex and has a reduced expression begining
with $w_{st}$ for some $t \in S$ with $st \ne ts$;}
\item{\rm (iv)}{there exists $I = \{u, u'\} \subset S$ with $s \not\in I$,
$uu' \ne u'u$, $su = us$ and $su' = u's$ such that when we write 
$w = w_I w^I$, we have $\ell(w_I) > 1$;}
\item{\rm (v)}{there exists $I = \{u, u'\} \subset S$ with 
$uu' \ne u'u$ such that when we write 
$w = (^Iw)(_Iw)$, we have $\ell(_Iw) > 1$;}
\item{\rm (vi)}{$x$ is weakly complex and there exist $t, u \in S$ with
$st \ne ts$, $ut \ne tu$ and $su = us$ such that $w$ has a reduced
expression of the form $$
u (tsts\cdots) x'
,$$ where the alternating product of $t$ and $s$ contains $m(s, t) - 1$ terms,
and we have $u(tuw) > tuw$;}
\item{\rm (vii)}{$x$ is weakly complex and there exist $t, u \in S$ with
$m(s, t) = 3$, $ut \ne tu$ and $su = us$ such that $w = sx$ has a reduced
expression of the form $w = utsux'$.}
\endproclaim

\demo{Proof}
This is \cite{{\bf 12}, Lemma 6.9}.
\qed\enddemo

\proclaim{Lemma \secc.8}
Let $W$ be a star reducible Coxeter group and let $x \in W$ be a
fully commutative or weakly complex element.  Then we have:
\item{\rm (i)}
{$\te_x \in \L$;}
\item{\rm (ii)}
{if $s \in S$ is such that $sx < x$, then $\te_x \in \latl{s}$;}
\item{\rm (iii)}
{if $s \in S$ is such that $xs < x$, then $\te_x \in \latr{s}$.}
\endproclaim

\demo{Proof}
The proof is by induction on $\ell(x)$, and the base case, $\ell(x) = 0$,
is easy.  In the inductive step, we will 
freely use the facts that, by Lemma \secc.6, the elements $sx$ and $xs$ 
occurring in assertions (ii) and (iii) satisfy the inductive hypotheses.

We first prove assertion (i).

If $\ell(x) > 0$, we may use a case analysis based on Lemma \secc.7 to prove
the first assertion.
If we are in case (i) of Lemma 3.7, this follows from the 
observation that if $x = s_1 s_2 \cdots s_r$ is a product of commuting 
generators, then $\te_x \in \latl{s}$ for each $s \in \{s_1, s_2, \ldots, 
s_r\}$.

In case (ii) of Lemma 3.7, we may assume that $x$ has a reduced expression 
beginning
with $st$, where $s$ and $t$ are noncommuting generators.  Since $tsx < sx$,
we have $\te_{sx} \in \latl{t}$ by induction, and then $\te_x \in \L$ by
Lemma \secc.5.  The analysis of case (iii) uses a similar argument.

In case (iv), we may assume that both $sx$ and $x$ have reduced expressions
beginning $uu'$, following the notation of Lemma \secc.7.  By induction,
$\te_{ux} \in \latl{u'}$, and hence $\te_x \in \L$ by Lemma \secc.5.
The analysis of case (v) uses a similar argument.

In case (vi), we have $x = uw_{st}x'$ reduced, so that $x$ has a reduced
expression beginning $ut$.  By induction, $\te_{ux} \in \latl{t}$, and 
hence $\te_x \in \L$ by Lemma \secc.5.  The analysis of case (vii) is the 
same, thus completing the proof of assertion (i).

We will now prove assertion (ii); the proof of assertion (iii) is by an 
analogous argument.

We know that $\te_{sx} \in \L$ by induction, and we have just shown that
$\te_x \in \L$.  Now $$
\te_x = \te_s \te_{sx} = (b_s - v^{-1}) \te_{sx}
,$$ and we have $v^{-1} \te_{sx} \in v^{-1} \L$ from the definitions,
which shows that $$
b_s \te_{sx} \in \L
.$$  By Lemma \secc.1, we have $$
b_s \te_{sx} = \sum_{w \in W_c} \l_w b_w
,$$ where $\l_w \ne 0$ implies $sw < w$, and the fact that 
$b_s \te_{sx} \in \L$ means that all $\l_w$ lie in $\A^-$.
Lemma \secc.4 shows that $b_s \te_{sx}$, and therefore $\te_s \te_{sx}$,
lies in $\latl{s}$, as required.
\qed\enddemo

\proclaim{Proposition \secc.9}
Let $W$ be a star reducible Coxeter group, let $s, t \in S$ be
noncommuting generators and let $w \in W_c$.   Then we have:
\item{\rm (i)}{$$
\te_s \te_w \in \cases
v \latl{s} & \text{ if } sw < w, \cr
\latl{s} & \text{ if } sw > w;
\endcases
$$} \item{\rm (ii)}
{$\te_s \L \cap \L \subseteq \latl{s}$;}
\item{\rm (iii)}
{$\te_s \latl{t} \subseteq \latl{st}$.}
\item{\rm (iv)}
{if $a \in S$ does not commute with $t$ and $a \ne s$, then 
$\te_a \latl{st} \subseteq \latl{a}$.}
\qed\endproclaim

\demo{Proof}
This was proved in \cite{{\bf 12}, Proposition 4.10} for any Coxeter group
satisfying the property that $\te_x \in \latl{u}$ whenever $x = uw$ is
a weakly complex element, $w \in W_c$ and $u \in S$.  This hypothesis
is satisfied by Lemma \secc.8 (ii).
\qed\enddemo

% 4. CF factorisation & reduced expressions

\head \secd. Main results \endhead

In \S\secd, we will show that any element of a star reducible Coxeter
group (not just a fully commutative element) has a reduced expression
with a particularly nice form.  More precisely, we have the following

\proclaim{Theorem \secd.1}
Let $W$ be a star reducible Coxeter group, and let $w \in W$.  Then
one of the following possibilities occurs for some Coxeter generators
$s, t, u$ with $m(s, t) \ne 2$, $m(t, u) \ne 2$ and $m(s, u) = 2$:
\item{\rm (i)}
{$w$ is a product of commuting generators;}
\item{\rm (ii)}
{$w$ has a reduced expression beginning with $st$;}
\item{\rm (iii)}
{$w$ has a reduced expression ending in $ts$;}
\item{\rm (iv)}
{$w$ has a reduced expression beginning with $sut$.}
%\item{\rm (v)}
%{$w$ has a reduced expression ending in $tsu$.}
\endproclaim

\demo{Proof}
Let $\bs$ be any reduced expression for $w$, and let $$
\bs_1 \bs_2 \cdots \bs_p
$$ be its Cartier--Foata normal form.  If $p = 1$, then case (i) applies,
and we are done.  

If not, let $t$ be a generator occurring in the factor
$\bs_2$.  By definition of the normal form, $t$ fails to commute with some
generator in $\bs_1$.  If $t$ fails to commute with only one such 
generator, $s$, then $\bs$ is commutation equivalent to a reduced expression
beginning with $st$, and case (ii) applies.

If $t$ fails to commute with precisely two generators, $s$ and $u$, in $\bs_1$,
then $\bs$ is commutation equivalent to a reduced expression beginning
$sut$, and we necessarily have $su = us$ by definition of the normal form,
so case (iv) applies.

Note that there cannot be four distinct generators $u_1, u_2, u_3, u_4$ in
$\bs_1$ not commuting with $t$, or $u_1 u_2 t u_3 u_4$ would be an element
of $W_c$ that is neither star reducible nor a product of commuting generators,
a contradiction.
We may therefore assume that each generator $t_i$ in $\bs_2$
fails to commute with precisely three (necessarily distinct and mutually 
commuting) generators, $\{u_{ij} : 1 \leq j \leq 3\}$, 
in $\bs_1$.  

Suppose that $\bs_2$ contains $k$ generators and the set $$
\{u_{ij} : 1 \leq i \leq k, \ 1 \leq j \leq 3 \}
$$ consists of $3k$ distinct elements of $\bs_1$.
This implies that, given such a $u_{ij}$, the only generator in $\bs_2$
not commuting with $u_{ij}$ is $t_i$.
Consequently, if $\bs_3$ is empty, then $w$ has a reduced
expression ending in $u_{11} t_1$, and case (iii) applies.  We may therefore
assume that $\bs_3$ contains a generator, $t'$.  We know $t'$ fails to commute
with some element of $\bs_2$, and without loss of generality, we may assume
that $m(t', t_1) \ne 2$.  None of the elements $\{t', u_{11}, u_{12}, u_{13}\}$
commutes with $t_1$, and if they were all distinct then $$
u_{11} u_{12} t_1 u_{13} t'
$$ would be an element of $W_c$ that would be neither star reducible nor
a product of commuting generators, a contradiction.  Without loss of 
generality, we may assume that $t' = u_{11}$, meaning that $w$ has a
reduced expression beginning $$
u_{12} u_{13} u_{11} t_1 u_{11}
.$$  If $m(t_1, u_{11}) = 3$, we may apply a braid relation to transform
this expression to one beginning $$
u_{12} u_{13} t_1 u_{11}
,$$ and case (iv) applies.  If, on the other hand, $m(t_1, u_{11}) > 3$,
the element $$
y = u_{12} u_{11} t_1 u_{11} u_{13}
$$ satisfies $y \in W_c$, but $y$ is neither star reducible nor
a product of commuting generators, a contradiction.

We have now reduced to the case where the set $$
\{u_{ij} : 1 \leq i \leq k, \ 1 \leq j \leq 3 \}
$$ is redundantly described.  Without loss of generality, we may assume
that $u := u_{11} = u_{21}$.  Now $$
y' = u_{12} u_{13} t_1 u t_2 u_{22} u_{23}
$$ lies in $W_c$, even if the set $\{u_{12}, u_{13}, u_{22}, u_{23}\}$
is redundantly described, because any two repeated occurrences of a generator
$s$ in the given reduced expression are separated by at least two occurrences
of generators not commuting with $s$ (see \cite{{\bf 11}, Remark 3.3.2}).
However, $y'$ is neither a product of commuting 
generators, nor star reducible, so this case cannot occur, completing the
analysis.
\qed\enddemo

\remark{Remark \secd.2}
By symmetry of the definitions, one can state a version of Theorem \secd.1
in which condition (iv) is replaced by the condition ``$w$ has a reduced 
expression ending in $tsu$''.
\endremark

The following result, which was proved by Losonczy \cite{{\bf 18}, Proposition 2.6,
Theorem 3.4} in type $D_n$, is new in type $E_n$, type $F_n$ ($n > 4$), type
$H_n$ ($n > 4$) and the two exceptional cases $\ti{E}_6$ and $\ti{F}_5$
discussed later (see Theorem \secf.3).

\proclaim{Theorem \secd.3}
If $W$ is a star reducible Coxeter group and 
$\L_\H$ is the free $\A^-$-submodule of $\H$ with basis 
$\{\T_w : w \in W\}$,
then the homomorphism $$\th : \H \ra TL(X)$$ restricts to an $\A^-$-linear map 
from $\L_\H$ to $\L$.
In particular, for any $w \in W$, we have $\th(\T_w) \in \L$, and 
$\pi(\th(\T_w)) = \pi(\th(C'_w))$.
\endproclaim

\demo{Proof}
We first prove that $\te_w \in \L$ using induction on $\ell(w)$ and the case
analysis of Theorem \secd.1.

If $w$ is a product of commuting generators, then $w \in W_c$ and the assertion
is immediate from the definitions.  This deals with the cases $\ell(w) \leq 1$.

If $w$ has a reduced expression beginning with $st$, as in Theorem \secd.1
(ii), then $\te_{sw}, \te_{tsw} \in \L$ by induction, and thus $$
\te_t \te_{tsw} = \te_{sw} \in \te_t \L \cap \L \subseteq \latl{t}
$$ by Proposition \secc.9 (ii).  We therefore have $$
\te_s \te_{sw} = \te_w \in \te_s \latl{t} \subset \L
$$ by Proposition \secc.9 (iii), as required.

If $w$ has a reduced expression ending in $ts$, as in Theorem \secd.1 (iii),
a symmetrical argument gives the desired conclusion.

Finally, suppose that $w$ has a reduced expression beginning with $sut$, 
as in Theorem \secd.1 (iv).  By induction, $\te_{uw}, \te_{suw}, \te_{tsuw}
\in \L$.  We also have $$
\te_t \te_{tsuw} = \te_{suw} \in \te_t \L \cap \L \subseteq \latl{t}
$$ by Proposition \secc.9 (ii), and $$
\te_s \te_{suw} = \te_{uw} \in \te_s \latl{t} \subset \latl{st}
$$ by Proposition \secc.9 (iii).  Finally, we have $$
\te_u \te_{uw} = \te_w \in \latl{u} \subset \L
$$ by Proposition \secc.9 (iv), as required.  

This completes the proof that $\te_w \in \L$, and it is then clear that
$\th(\T_w) \in \L$.  Since $C'_w$ and $\T_w$ agree modulo $v^{-1} \L_\H$
(as explained in, for example, \cite{{\bf 14}, Proposition 1.2.2}), the final
claim also follows.
\qed\enddemo

\proclaim{Lemma \secd.4}
Let $W$ be an arbitrary Coxeter group such that $I = \{s, t\} 
\subseteq S$ is a pair of noncommuting generators, and suppose that 
$w \in W$ satisfies $tw < w$ and $sw > w$.
\item{\rm (i)}
{If $w \in W_c$, $tw < w$ and $sw \not\in W_c$, then $sw = w_{st}w'$ is 
reduced.}
\item{\rm (ii)}
{Taking star operations with respect to $I$, we have $$
C'_s C'_w = C'_{^*w} + C'_{_*w} \mod J(X)
,$$ where $C'_z$ is defined to be zero if $z$ is an undefined symbol.}
\endproclaim

\demo{Note}
There is also a right-handed version of this result.
\enddemo

\demo{Proof}
Part (i) follows from \cite{{\bf 19}, Proposition 2.3} (see also \cite{{\bf 12}, Lemma
4.5 (i)}), and part (ii) follows from \cite{{\bf 12}, Lemma 6.2}.
\qed\enddemo

\proclaim{Lemma \secd.5}
Let $W$ be a star reducible Coxeter group and let $x \in W$ be weakly
complex.  Then $\th(C'_x) = 0$.
\endproclaim

\demo{Proof}
We write $x = sw$ with $s \in S$ and $w \in W_c$.  The proof is by 
induction on $\ell(x)$, using Lemma \secc.6 and the case analysis of Lemma 
\secc.7.

Since $x$ is weakly complex, we are in one of cases (iii)--(vii) of Lemma
\secc.7.  Let us first suppose we are in case (iii), meaning that 
$x = w_{st} x'$ is reduced.  Since $x$ has a reduced expression beginning
with $st$ and $\te_x, \te_{sx} \in \L$ by Theorem \secd.3, Proposition
\secc.9 (ii) shows that $\te_x \in \latl{s}$.  Similarly, $x$ has a reduced
expression beginning with $ts$, and $\te_x \in \latl{t}$.  Since $s$ and $t$
do not commute, we have $\latl{s} \cap \latl{t} \subseteq v^{-1} \L$,
which shows that $\pi(\te_x) = 0$.  By Theorem \secd.3, we also have 
$\pi(\th(C'_x)) = 0$.  Since $\overline{\th(C'_x)} = \th(C'_x)$ and $\th(C'_x)
\in \L$, \cite{{\bf 15}, Lemma 2.2.2} shows that $\th(C'_x) = 0$, as required.

Suppose that we are in case (iv) of Lemma \secc.7.
We may assume without loss of generality that $w$ and $x$ each have a reduced
expression beginning $uu'$, where $I' = \{u, u'\}$ is a pair of noncommuting
generators and $u, u'$ satisfy the conditions of Lemma \secc.7 (iv).  
We cannot have $x = w_{uu'} x'$ reduced, or $w = sx = w_{uu'} (sx')$ would not
be fully commutative, which is a contradiction.
Taking star operations with respect to $I'$, we may therefore assume that
${^*ux}$ is defined and equal to $x$, and furthermore (by Lemma \secd.4 (i)), 
that $ux$ is weakly complex.
By Lemma \secd.4 (ii), we then have $$
C'_u C'_{ux} = C'_x + C'_{_*ux} \mod J(X)
.$$  Since $C'_{ux} \in J(X)$ by induction, we need to show that $C'_{_*ux}
\in J(X)$.  We may assume that ${_*ux}$ is defined, or this is obvious.
By Lemma \secc.6, either ${_*ux}$ is weakly complex or fully commutative,
and in the former case we are done by the inductive hypothesis.  However,
if ${_*ux} \in W_c$, then the fact that $ux \not\in W_c$ implies by
Lemma \secd.4 (i) that $u.ux = x < ux$, a contradiction.  This completes
the analysis of case (iv), and case (v) follows by a similar argument.
The only difference in the argument needed to treat case (v) is that we may
have $x = x' w_{uu'}$ reduced, in which case we are done by an argument like
that used to treat case (iii).

Suppose we are in case (vi) of Lemma \secc.7, and keep the same notation.
In this case, we have $x = uw_{st}x'$ reduced, and furthermore, $w_{st}$
has a reduced expression beginning with $t$, which does not commute with $u$.
As in case (iii), we may assume that we do not have $x = w_{tu}x'$ reduced.
Taking star operations with respect to $I'' = \{u, t\}$,
we may assume as in the analysis of case (iv) that ${^*ux}$ is defined and 
equal to $x$, and that $ux$ is weakly complex.
By Lemma \secd.4 (ii), we now have $$
C'_u C'_{ux} = C'_x + C'_{_*ux} \mod J(X)
.$$  As in the analysis of case (iv), the only nonobvious case left 
to consider is when ${_*ux}$ is defined and fully commutative.  In this 
case, ${_*ux}$ is reduced of the form $$
(stst\cdots)x'
,$$ where there are $m(s, t) - 1$ occurrences of $s$ or $t$.  However, this
cannot happen: $t \in \ldescent{ux}$ implies that $u \in \ldescent{_*ux}$,
and a fully commutative element cannot have a reduced expression beginning
with $st$ and another beginning with $u$ if $s \ne u$ and $t$ and $u$ do 
not commute.

The analysis for case (vii) is exactly the same as that for case (vi), and
this completes the proof.
\qed\enddemo

\proclaim{Theorem \secd.6}
Let $W$ be a star reducible Coxeter group.
\item{\rm (i)}
{If $w \in W$ is weakly complex, then $\te_w \in v^{-1} \L$; in other
words, $W$ has ``Property W'', in the sense of \cite{{\bf 12}}.}
\item{\rm (ii)}
{If $w \in W_c$, then we have $$
c_s c_w = \cases
(v + v^{-1}) c_w & \text{ if } \ell(sw) < \ell(w),\cr
c_{sw} + \sum_{sy < y} \mu(y, w) c_y
& \text{ if } \ell(sw) > \ell(w),\cr
\endcases$$ where $c_z$ is defined to be zero whenever $z \not\in W_c$, and
where $\mu(y, w)$ is the integer defined by Kazhdan and Lusztig 
in \cite{{\bf 17}}.}
\item{\rm (iii)}
{If $I = \{s, t\}$ is a pair of noncommuting generators,
and we have $w \in W_c$ with $tw < w$, then we have $$
c_s c_w = c_{^*w} + c_{_*w}
,$$  where $c_z$ is defined to be zero whenever $z$ is an undefined 
symbol.}\item{\rm (iv)}
{If $w \in W_c$, then $c_w = \th(C'_w)$.}
\item{\rm (v)}
{The structure constants arising from the $c$-basis of $TL(X)$ lie in 
$\zed^{\geq 0}[\d]$.}
\qed\endproclaim

\demo{Proof}
For part (i), 
let $w \in W$ be a weakly complex element.  We know from Theorem \secd.3 
that $\pi(\th(\T_w)) = \pi(\th(C'_w))$, and we know from Lemma \secd.5 that
$\pi(\th(C'_w)) = 0$.  Part (i) is immediate from these observations.

Part (ii) is essentially \cite{{\bf 12}, Theorem 5.13}, the only difference
being that (i) allows us to remove the extra hypothesis that $W$ should
have Property W.  Similarly, parts (iii) and (iv) now follow from \cite{{\bf 12}, 
Proposition 6.3}, and part (v) now follows from \cite{{\bf 12}, Theorem 6.13}.
\qed\enddemo

\remark{Remark \secd.7}
Note that part (iii) of the theorem allows the $c$-basis to be constructed
inductively.  Part (v) proves \cite{{\bf 14}, Conjecture 1.2.4} for star reducible
Coxeter groups.  This is a new result for type $F_n$ ($n > 4$) and type
$\ti{F}_5$ (see Lemma \sece.5), and it provides a new elementary
proof of positivity in type $\ti{C}_{n-1}$ (for $n$ even).
\endremark

% 5. Examples

\head \sece. Some examples of star reducible Coxeter groups \endhead

In \S\sece, we present some specific examples of star reducible Coxeter groups,
and we present various methods to construct new examples out of known ones.
It will turn out in \S\secf\  that these methods suffice to construct all
examples, assuming as always that the Coxeter generating set $S$ is finite.

In order to show that certain Coxeter groups are star reducible, we need to
associate a sequence of graphs to each Cartier--Foata normal form.  This idea
has also been used by Fan in \cite{{\bf 5}, Lemma 4.3.2}, and by Fan and the author
in \cite{{\bf 6}, \S2.4}.

\definition{Definition \sece.1}
Let $\bs$ be an element of the commutation monoid $\coms$ with Cartier--Foata
normal form $\bs = \bs_1 \bs_2 \cdots \bs_p
.$  For all $1 \leq i < p$, we define the graph $X_i(\bs)$ to be the induced
labelled subgraph of $X$ corresponding to the set of all generators appearing 
in the factors $\bs_i$ and $\bs_{i+1}$.  If $w \in W_c$, then we define
$X_i(w)$ to be the graph $X_i(\bs)$, where $\bs$ is the (unique) element 
of $\coms$ corresponding to $w$.
\enddefinition

\remark{Remark \sece.2}
If $\bs$ is a reduced expression for some Coxeter group element, 
the generators appearing in the subword $\bs_i \bs_{i+1}$ of $\bs$ 
are distinct, by definition of the normal form.
\endremark

For the next lemma, we assume that the Coxeter group $(W, S)$ is of type
$\ti{C}_{\twola}$ ($l \geq 1$), meaning that 
$S = \{s_1, s_2, \ldots, s_{\twolb}\}$ and we have the relations

\item{\rm (a)}
{$m(s_i, s_j) = 2$ if $|i - j| > 1$,}
\item{\rm (b)}
{$m(s_1, s_2) =  m(s_{\twola}, s_{\twolb}) = 4$,}
\item{\rm (c)}
{$m(s_i, s_{i+1}) = 3$ if $1 < i < \twola$.}

\proclaim{Lemma \sece.3}
Let $W$ be the Coxeter group of type $\ti{C}_{\twola}$, with the
above notation.
Suppose that $\bs \in \coms$ corresponds to a reduced expression for 
$w \in W_c$, and let $\bs_1 \bs_2 \cdots \bs_p$ be the Cartier--Foata
normal form of $\bs$.  Suppose also that $w \in W_c$ is not left star 
reducible.
Then, for $1 \leq i < p$ and $1 \leq j \leq \twolb$, the following hold:
\item{\rm (i)}
{if $s_1$ occurs in $\bs_{i+1}$, then $s_2$ occurs in $\bs_i$;}
\item{\rm (ii)}
{if $s_{\twolb}$ occurs in $\bs_{i+1}$, then $s_{\twola}$ occurs in $\bs_i$;}
\item{\rm (iii)}
{if $j \not\in \{1, \twolb\}$ and $s_j$ occurs in $\bs_{i+1}$, then both 
$s_{j-1}$ and $s_{j+1}$ occur in $\bs_i$.}
\endproclaim

\demo{Proof}
The assertions of (i) and (ii) are immediate from properties of the normal
form, because $s_2$ (respectively, $s_{\twola}$) is the only generator not
commuting with $s_1$ (respectively, $s_{\twolb}$).  We will now prove (iii) 
by induction on $i$.  Suppose first that $i = 1$.

Suppose that $j \not\in \{1, \twolb\}$ and that $s_j$ occurs in $\bs_2$.  
By definition of the normal form, there must be a generator $s \in \bs_1$ 
not commuting with $s_j$.   Now $s$ cannot be the only such generator, or 
$w$ would be left star reducible to $sw < w$. 
Since the only generators not
commuting with $s_j$ are $s_{j-1}$ and $s_{j+1}$, these must both
occur in $\bs_1$.

Suppose now that the statement is known to be true for $i \leq N$, and
let $i = N+1 \geq 2$.
Suppose also that $j \not\in \{1, \twolb\}$ and $s_j$ occurs in $\bs_{N+1}$.  
As in the base case, there must be at least one generator $s$ occurring 
in $\bs_N$ that does not commute with $s_j$.  

Let us first consider the case where $j \not\in \{2, \twola\}$, and write
$s = s_k$ for some $1 \leq k \leq \twolb$.  The restrictions on $j$ mean
that $2 < k < \twola$.  By the inductive hypothesis,
this means that $s_{k-1}$ and $s_{k+1}$ both occur in
$\bs_{N-1}$, and that $m(s_{k-1}, s_k) = m(s_k, s_{k+1}) = 3$.
Now either $j = k-1$ or $j = k+1$; we consider the first possibility,
the other being similar.  (Since $j \geq 3$, this means $k \geq 4$.)  If
$s_{k-2}$ occurs in $\bs_N$, then statement (i) follows as both generators
not commuting with $s_j$ lie in $\bs_N$.  If, on the other hand, $s_{k-2}$
does not occur in $\bs_N$, the fact that $s_{k-1}$ occurs both in $\bs_{N-1}$
and in $\bs_{N+1}$ means that the word $\bs$ can be parsed in the form $
\bu_1 s_{k-1} \bu_2 s_{k-1} \bu_3
,$ where all the generators in $\bu_2$ commute with $s_{k-1}$ except for
one occurrence of $s_k$.  This means that $\bs$ is represented by a word
in $S^*$ containing a subword $s_{k-1} s_k s_{k-1}$, which contradicts
the assumption $w \in W_c$.

Now suppose that $j = 2$ (the case $j = \twola$ follows by a symmetrical
argument).  If both $s_1$ and $s_3$ occur in $\bs_N$, then
we are done.  If $s_3$ occurs in $\bs_N$ but $s_1$ does not, then the
argument of the previous paragraph applies.  Suppose then that $s_1$ occurs
in $\bs_N$ but $s_3$ does not.  By statement (i), $s_2$ occurs in 
$\bs_{N-1}$, and we cannot have $N = 2$,
or $w$ would be left star reducible to $s_2 w < w$.  Applying the inductive
hypothesis to (i), we see that $s_1$ and $s_3$ both occur in $\bs_{N-2}$.
Putting all this together, we find that $w$ has a reduced expression
containing a subword of the form $s_3 s_1 s_2 s_1 s_2$, which is incompatible
with $w \in W_c$.  This completes the inductive step.
\qed\enddemo

\proclaim{Proposition \sece.4}
A Coxeter group of type $\ti{C}_{\twola}$ is star reducible.
\endproclaim

\demo{Proof}
Keeping the previous notation, we suppose that $w \in W_c$ is not left star 
reducible and prove that either $w$ is a product of commuting generators, 
or $w$ is right star reducible.

If $\bs_2$ is empty, then $w$ is a product of commuting generators, and
we are done.  Otherwise, the graph $X_{p-1}(w)$ has the property that not
all of its connected components have size $1$.  Let $\Gamma$ be one of
the components with $|\Gamma| > 1$.

Suppose that $\Gamma = \ti{C}_{\twola}$, which has an even
number of vertices.  Either this forces $s_1$ to occur in $\bs_{p-1}$ and
$s_2$ to occur in $\bs_p$, or it forces $s_{\twolb}$ to occur in $\bs_{p-1}$
and $s_{\twola}$ to occur in $\bs_p$.  In the first case, $w$ is right star
reducible with respect to $\{s_1, s_2\}$, and in the second, $w$ is right 
star reducible with respect to $\{s_{\twola}, s_{\twolb}\}$.

Suppose that $\Gamma$ is a Coxeter graph of type $B_n$.  Conditions (i)--(iii)
of Lemma \sece.3 show that there are four possibilities:
\item{\rm (a)}
{$s_1$ occurs in $\bs_{p-1}$ and corresponds to a vertex of $\Gamma$, 
and $n$ is odd; }
\item{\rm (b)}
{$s_{\twolb}$ occurs in $\bs_{p-1}$ and corresponds to a vertex of $\Gamma$, 
and $n$ is odd; }
\item{\rm (c)}
{$s_1$ occurs in $\bs_p$ and corresponds to a vertex of $\Gamma$, 
and $n$ is even; }
\item{\rm (d)}
{$s_{\twolb}$ occurs in $\bs_p$ and corresponds to a vertex of $\Gamma$, 
and $n$ is even.}

Let $k > 1$ be the number of vertices in $\Gamma$.  In case (a), $w$ is
right star reducible with respect to $\{s_1, s_2\}$, and in case (b), with
respect to $\{s_{\twola}, s_{\twolb}\}$.  In case (c), $w$ is right star 
reducible with respect to $\{s_{k-1}, s_k\}$, and in case (d), with respect to
$\{s_{\twolc-k}, s_{\twold-k}\}$.

The only other possibility is that $\Gamma$ is a Coxeter graph of type
$A_k$.  In this case, condition (iii) of Lemma \sece.3 forces $k > 1$ to
be odd.  If $\{s_a, s_{a+1}, \ldots, s_b\}$ are the generators involved
in $\Gamma$, then $s_a$ and $s_b$ both lie in $\bs_{p-1}$, and $w$ is right 
star reducible with respect to $\{s_a, s_{a+1}\}$ and with respect to 
$\{s_{b-1}, s_b\}$.
\qed\enddemo

\proclaim{Lemma \sece.5}
Let $W$ be the Coxeter group with Coxeter matrix $$
(m_{i, j})_{1 \leq i, j \leq 6}= \left( \matrix
1 & 3 & 2 & 2 & 2 & 2 \cr
3 & 1 & 3 & 2 & 2 & 2 \cr
2 & 3 & 1 & 4 & 2 & 2 \cr
2 & 2 & 4 & 1 & 3 & 2 \cr
2 & 2 & 2 & 3 & 1 & 3 \cr
2 & 2 & 2 & 2 & 3 & 1 \cr
\endmatrix \right)
,$$ and denote $S = \{s_1, s_2, \ldots, s_6\}$ in the obvious way.
Then $W$ is star reducible.
\endproclaim

\topinsert
\topcaption{Figure 1} The Coxeter graph $X$ of Lemma \sece.5
\endcaption
\centerline{
\hbox to 3.069in{
\vbox to 0.305in{\vfill
        \includegraphics{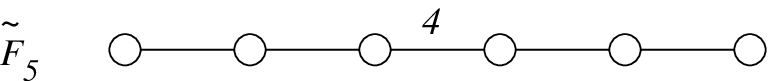}
}
\hfill}
}
\endinsert

\demo{Note}
The graph $X$ in this case is shown in Figure 1.
Note that there is a symmetry of the graph
$X$, namely that sending $s_i$ to $s_{7-i}$, which induces a Coxeter group
automorphism of $W(X)$.
\enddemo

\demo{Proof}
Let $w \in W_c$ be such that $w$ is not left star reducible or a product of 
commuting generators, and suppose (for a contradiction) that $w$ is not 
right star reducible.

Let $\bs \in \coms$ correspond to $w$, and let $$
\bs = \bs_1 \bs_2 \cdots \bs_p
$$ be the corresponding Cartier--Foata normal form.  Since $w$ is not a product
of commuting generators, there exists a generator $s_k \in \bs_2$.  Since $w$
is not left star reducible, there must be at least two generators in
$\bs_1$ that do not commute with $s$.  Because $X$ is a straight line, these
two generators must be $s_{k-1}$ and $s_{k+1}$, so that in particular we
cannot have $k = 1$ or $k = 6$.  Since $|X| = 6$ and the generators from 
$\bs_1$ pairwise commute, we must therefore have $2 \leq |\bs_1| \leq 3$.

Suppose first that $|\bs_1| = 2$.  By symmetry of $X$ and the above remarks,
it suffices to consider the cases $\bs_1 = s_1 s_3$ and $\bs_1 = s_2 s_4$.

If $\bs_1 = s_1 s_3$ then $\bs_2$ can only contain $s_2$, for 
if $\bs_2$ contained $s_4$ (the only other generator not commuting with 
either $s_1$ or $s_3$) then $w$ would be 
left star reducible to $s_3 w$.  Now $s_1 s_3 s_2$ is right star reducible,
so $\bs_3$ must contain a generator, and this generator must not commute with
$s_2$.  We cannot have $s_1$ occurring in $\bs_3$, or $w$ would have a 
reduced expression containing $s_1 s_2 s_1$ consecutively.
Similarly, we cannot have $s_3$ occurring in $\bs_3$, producing a 
contradiction.

If $\bs_1 = s_2 s_4$ then, arguing as in the above paragraph, we find that
$\bs_2 = s_3$, $\bs_3 = s_4$, $\bs_4 = s_5$ and $\bs_5 = s_6$.  At this point,
we are stuck, and $s_2 s_4 s_3 s_4 s_5 s_6$ is right star reducible, which
is a contradiction.

Suppose now that $|\bs_1| = 3$.  By symmetry of $X$, we may assume that
$\bs_1 = s_1 s_3 s_5$.  Now $\bs_2$ is nonempty, but it cannot contain $s_6$,
or $w$ would be left star reducible to $s_5 w$.  If $\bs_2$ contains only
$s_4$, then $\bs_3$ must be nonempty as $s_1 s_3 s_5 s_4$ is right star
reducible.  In turn, we must have $\bs_3 = s_3$, $\bs_4 = s_2$,
$\bs_5 = s_1$, and then there are no possible choices
for $\bs_6$, a contradiction.  If $\bs_2$ contains 
only $s_2$, then a similar argument shows that all choices for $\bs_3$ lead
to a contradiction.  The only other possibility is for $\bs_2 = s_2 s_4$,
which forces $\bs_3 = s_3$.  However, $s_1 s_3 s_5 s_2 s_4 s_3$ is right
star reducible, and we must then have $\bs_4 = s_4$, $\bs_5 = s_5$ and
then there are no possible choices for $\bs_6$, a contradiction.
We have exhausted all the possibilities, so the
assumption that $w$ is not right star reducible is wrong, completing the 
proof.
\qed\enddemo

% FC finite

\proclaim{Lemma \sece.6}
If $W$ is a Coxeter group for which $W_c$ is finite, then $W$ is
star reducible.
\endproclaim

\demo{Proof}
As pointed out in \cite{{\bf 12}, Remark 3.5}, this result follows from the
argument of \cite{{\bf 5}, Lemma 4.3.1} together with \cite{{\bf 19}, Proposition 2.3}.

The way the argument works is as follows.  Suppose that $w_1 \in W_c$ has
the property that $w_1$ is neither left nor right star reducible.  Let
$\bs$ be a reduced expression for $w_1$, and let $$
\bs = \bs_1 \bs_2 \cdots \bs_p
$$ be the Cartier--Foata normal form of the corresponding element of
$\coms$.  The results of Fan and Stembridge just mentioned show that $$
\bs_p \bs_{p-1} \cdots \bs_2 \bs_1 \bs_2 \cdots \bs_{p-1} \bs_p
$$ is also a reduced expression for an element $w_2 \in W_c$ that also
has the property that it cannot be left or right star reduced.  Proceeding
in this way, we obtain an infinite sequence $\{w_i\}_{i \in \enn}$ of
distinct elements of $W_c$, which contradicts the hypothesis.
\qed\enddemo

% complete graph

\proclaim{Lemma \sece.7}
If $W$ is a Coxeter group for which no two distinct elements of $S$
commute, then $W$ is star reducible.
\endproclaim

\demo{Proof}
Let $w \in W_c$.  The hypotheses show that $w$ has a unique reduced 
expression, $$
w = s_1 s_2 \cdots s_k
.$$  Since $s_1$ and $s_2$ do not commute by hypothesis, $w$ is star
reducible to $s_1 w < w$.  Iterating this argument proves the assertion.
\qed\enddemo

% parabolic/connected component

The following useful lemma is an easy consequence of the definitions.

\proclaim{Lemma \sece.8}
If $(W, S)$ is a star reducible Coxeter group, then so is any parabolic
subgroup $(W_I, I)$ of $(W, S)$.  In particular, any connected component
of the Coxeter graph of a star reducible Coxeter group corresponds to another
star reducible Coxeter group.
\qed\endproclaim

% reduction

\definition{Definition \sece.9}
Let $(W, S)$ be a Coxeter group corresponding to Coxeter graph $X$ and
function $m : S \times S \ra \enn$.
We define the Coxeter group $(\upsilon(W), S)$ to be the group corresponding
to the function $m' : S \times S \ra \enn$, where $$
m'(i, j) = \cases
m(i, j) & \text{ if } m(i, j) < 3;\cr
3 & \text { otherwise.}\cr
\endcases
$$  In other words, it is the group obtained by deleting all edge labels
bigger than $3$ (including edges with infinite label) in $X$.
\enddefinition

\proclaim{Lemma \sece.10}
If $(W, S)$ is a star reducible Coxeter group, then so is $(\upsilon(W), S)$.
\endproclaim

\demo{Proof}
Let $w \in \upsilon(W)$ be a fully commutative element, and
let $\bs$ be a reduced expression for $w$.  Since all reduced expressions
for $w$ are commutation equivalent, and since two generators $s, s' \in S$
commute in $\upsilon(W)$ if and only if they commute in $W$, it follows
that $\bs$ is also a reduced expression for a fully commutative element 
$w^+ \in W$.

Since $W$ is star reducible, either $w^+$ is a product of commuting generators
in $W$ (which means that $w$ is a product of commuting generators in 
$\upsilon(W)$), or $w^+$ is left or right star reducible to some other element 
of $W$.  We treat the case of left star reducibility, since the other case
is similar.  Suppose that $w^+$ is left star reducible with respect to 
$I = \{s, s'\} \subseteq S$.  If $m, m' : S \times S \rightarrow \enn$ are
the functions arising from the Coxeter groups $(W, S)$ and $(\upsilon(W), S)$
respectively, then we have $m(s, s') \geq m'(s, s') \geq 3$ by Definition
\sece.9 and the fact that $s, s'$ do not commute.  This means that we can
identify the $\{s, s'\}$-string, $S_w$, in $\upsilon(W)$ containing $w$ with a 
subset of the $\{s, s'\}$-string, $S_{w^+}$, in $W$ containing $w^+$; here
$S_w$ will consist of the $(m'(s, s') - 1)$ shortest elements of $S_{w^+}$.
Since star reducibility moves $w^+$ to a shorter element in $S_{w^+}$, there
is a corresponding star reduction of $w$ to a shorter element in $S_w$.
By iterating this procedure, we see that $w$ can be star reduced to a product
of commuting generators, as required.
\qed\enddemo

% list of the reduced case

The benefit of Lemma \sece.10 is that the simply laced star reducible Coxeter
groups have already been classified \cite{{\bf 11}}.

\proclaim{Theorem \sece.11 \cite{{\bf 11}}}
Let $W$ be a simply laced Coxeter group with (finite) generating set $S$.  
Then $W$ is star reducible if and only if each component of $X$ is either
a complete graph $K_n$ or appears
in the list depicted in Figure 2: type $A_n$ ($n \geq 1$), type $D_n$ 
($n \geq 4$), type $E_n$ ($n \geq 6$), type $\ti{A}_{n-1}$ ($n \geq 3$ and
$n$ odd) or type $\ti{E}_6$.
\endproclaim

\topinsert
\topcaption{Figure 2} Connected incomplete graphs associated to simply 
laced star reducible Coxeter groups
\endcaption
\centerline{
\hbox to 3.319in{
\vbox to 4.638in{\vfill
        \includegraphics{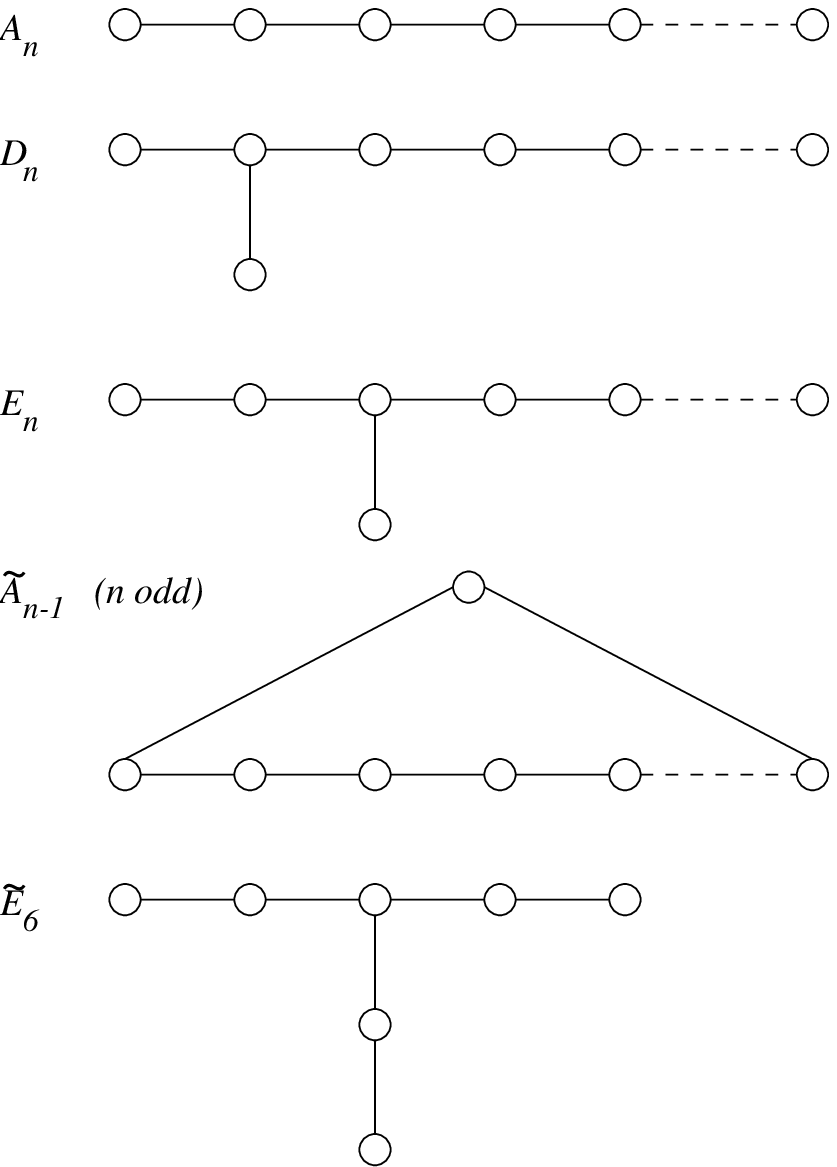}
}
\hfill}
}
\endinsert

\demo{Note}
The corresponding result for arbitrary $|S|$ is not much more difficult,
but we do not state it in order to avoid cardinality issues.
\enddemo

\demo{Proof}
This is a restatement of \cite{{\bf 11}, Theorem 1.5.2} using the definitions
and remarks of \cite{{\bf 11}, \S1.2}.
\qed\enddemo

% 6. Classification (partial reduction to simply laced case)

\head \secf. Classification of star reducible Coxeter groups \endhead

We are now ready to classify the star reducible Coxeter groups $(W, S)$
for finite $S$.  During the argument, which is reminiscent of the
classification of finite Coxeter groups \cite{{\bf 16}, \S2} and the classification
of FC-finite Coxeter groups (see \cite{{\bf 19}, \S4}, \cite{{\bf 8}, \S7}), we
will freely use the contrapositive statement to Lemma \sece.8.

By Lemma \sece.10 and Theorem \sece.11, the remaining part of this task
will be to determine how the edge labels in the graphs listed in Figure 2
may be increased so as to obtain another star reducible Coxeter group.
We first deal with the case where the graph has a branch point, which means
that it is of type $D_n$, $E_n$ or $\ti{E}_6$.

\proclaim{Lemma \secf.1}
Suppose that $X$ is a connected Coxeter graph with a branch point, and that
$W(X)$ is star reducible.  Then $X$ is simply laced.
\endproclaim

\topinsert
\topcaption{Figure 3} Coxeter graphs considered in the proof of Lemma 
\secf.1
\endcaption
\centerline{
\hbox to 2.888in{
\vbox to 0.777in{\vfill
        \includegraphics{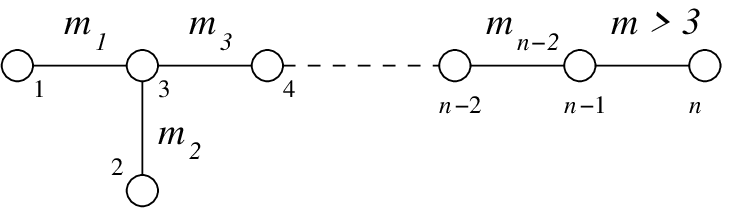}
}
\hfill}
}
\endinsert

\demo{Proof}
By the remarks preceding the statement (and Lemma \sece.8), it is enough to
show that $X$ cannot arise from a graph of Coxeter type $D_n$, where the label
of the edge furthest from the branch point is greater than $3$, and where
some of the other edges with labels $m \geq 3$ may also have been increased;
see Figure 3.
(If $n = 4$, the condition is that at least one of the edge labels must
strictly exceed $3$.)

Labelling the vertices as in Figure 3 (where vertices $1$ and $2$ commute,
$3$ is the branch point, and $m(n-1, n) > 3$), we find that $$
(s_1 s_2) s_3 s_4 \cdots s_{n-2} s_{n-1} s_n s_{n-1} s_{n-2} \cdots
s_4 s_3 (s_1 s_2)
$$ is a fully commutative element that cannot be left or right star reduced,
but that is not a product of commuting generators, which completes the
proof.
\qed\enddemo

\proclaim{Lemma \secf.2}
Suppose that $X$ is a Coxeter graph whose unlabelled graph is a $k$-cycle,
where $k \geq 5$ is odd, and that $W(X)$ is star reducible.  Then $X$ is 
simply laced.
\endproclaim

\demo{Proof}
Numbering the Coxeter generators $s_1, s_2, \ldots, s_k$ in an obvious
cyclic fashion, let us assume that $m(s_k, s_1) > 3$.  Since $k \geq 5$,
we have $m(s_2, s_k) = 2$ and $m(s_1, s_{k-1}) = 2$.  In this case, the
element $$
(s_2 s_k) s_1 s_k s_{k-1} s_{k-2} \cdots s_3 s_2 s_1 s_k (s_{k-1} s_1)
$$ is a fully commutative element that cannot be left or right star reduced,
but that is not a product of commuting generators.
\qed\enddemo

Finally, we may classify all star reducible Coxeter groups with a finite
generating set.

\topinsert
\topcaption{Figure 4} 
Connected incomplete graphs associated to non simply 
laced star reducible Coxeter groups
\endcaption
\centerline{
\hbox to 3.319in{
\vbox to 3.055in{\vfill
        \includegraphics{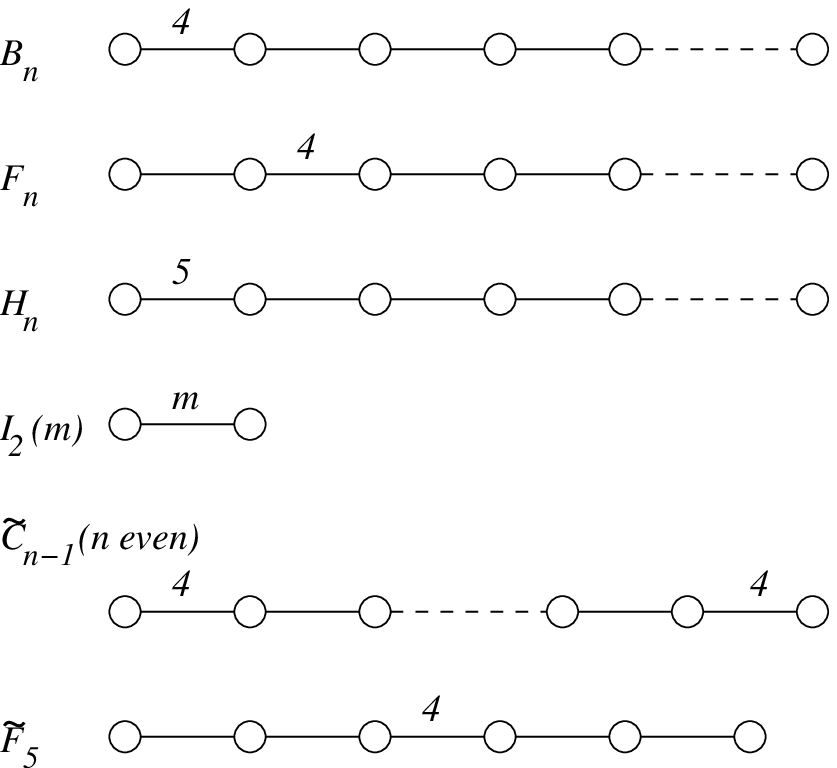}
}
\hfill}
}
\endinsert

\proclaim{Theorem \secf.3}
Let $W(X)$ be a Coxeter group with (finite) generating set $S$.  
Then $W(X)$ is star reducible if and only if each component of $X$ is either
a complete graph with all labels $m(i, j) \geq 3$, or appears
in one of the lists depicted in Figure 2 or Figure 4: 
type $A_n$ ($n \geq 1$), type $B_n$ ($n \geq 2$), type $D_n$ 
($n \geq 4$), type $E_n$ ($n \geq 6$), type $F_n$ ($n \geq 4$), type 
$H_n$ ($n \geq 2$), type $I_2(m)$ ($m \geq 3$), 
type $\ti{A}_{n-1}$ ($n \geq 3$ and
$n$ odd), type $\ti{C}_{n-1}$ ($n \geq 4$ and $n$ even), type $\ti{E}_6$
or type $\ti{F}_5$.
\endproclaim

\demo{Note}
Although there appear to be ten infinite families in the classification above,
the family $I_2(m)$ consists entirely of complete graphs and may thus be
incorporated into another family.
\enddemo

\demo{Proof}
We first summarize why the examples listed are star reducible.
The families $A$, $B$, $D$, $E$, $F$, $H$, $I$ have the property that $W_c$
is finite (see \cite{{\bf 19}, \S4}, \cite{{\bf 8}, \S7}), so they are star reducible
by Lemma \sece.6.  Types $\ti{A}_{n-1}$ and $\ti{E}_6$ are covered by 
Theorem \sece.11, type $\ti{C}_{n-1}$ is covered by Proposition \sece.4,
and type $\ti{F}_5$ is covered by Lemma \sece.5.

Let us now prove that the list given is complete, bearing in mind that 
Lemma \sece.8 allows us to reduce consideration to connected components.
If $W$ is star reducible, Lemma \sece.10 shows that $\upsilon(W)$ is as well.
If the graph $X$ is complete, then any increased labels are permissible by
Lemma \sece.7, so our list of complete graphs is correct.

There is no way to increase the labels of edges of the graphs of types $D$,
$E$ or $\ti{E}_6$ appearing in Figure 2 by Lemma \secf.1, so our list of
graphs with branch points is complete.  

If the Coxeter graph $X$ is a cycle
and $W$ is star reducible, it must be a cycle of odd length by Lemma 
\sece.10 and Theorem \sece.11.  A cycle of length $3$ is a complete graph, 
and then any labels are
permissible.  A cycle of length $5$ or greater cannot have any labels
increased by Lemma \secf.2, so our list of cycle shaped graphs is complete.

We have reduced consideration to the case where $X$ is a straight line.
Let us label the Coxeter generators $s_1, s_2, \ldots, s_n$ in an obvious
way.  We shall assume that $n \geq 3$, or else $X$ is complete, which we have
dealt with above.

We first show that $X$ has no edge labelled $6$ or greater.  To check this,
it is enough by Lemma \sece.8 to consider the case where $n = 3$ and
$m(s_2, s_3) \geq 6$.  In this case, the element $$
s_1 s_3 s_2 s_3 s_2 s_1 s_3
$$ provides the required counterexample of a fully commutative element that
is not a product of commuting generators, but also not left or right star
reducible.

Suppose now that $X$ has an edge labelled $5$ (but no labels strictly 
greater than $5$, by the above).  We claim that this edge must
be extremal.  If not, we may reduce to the case where $n = 4$ and 
$m(s_2, s_3) = 5$.  In this case, $$
s_1 s_3 s_2 s_3 s_2 s_4
$$ provides the required counterexample.

Suppose that $X$ has an extremal edge labelled $5$.  In this case, we claim
that this edge is the only edge with a label greater than $3$.  If not, we
may reduce (using Lemma \sece.8 as always) to the case where $m(s_1, s_2) = 5$
and $m(s_{n-1}, s_n) > 3$.  In this case, the element $$
s_1 s_3 s_2 s_1 s_2 s_3 s_4 \cdots s_{n-1} s_n s_{n-1} \cdots
s_4 s_3 s_2 s_1 s_2 s_1 s_3
$$ provides the required counterexample.  We conclude that if $X$ has an 
edge with label $5$, then $X$ is of type $H_n$, which is on the list.

Suppose now that $X$ has at least two edges labelled $4$, but no edge with 
label $5$ or higher.  If one of these edges is not extremal,
then we may reduce to the case where $m(s_2, s_3) = 4$ and 
$m(s_{n-1}, s_n) = 4$, and $$
s_1 s_3 s_2 s_3 s_4 \cdots s_{n-1} s_n s_{n-1} \cdots
s_4 s_3 s_2 s_1 s_3
$$ provides the required counterexample.  We deduce that there are precisely
two edges labelled $4$, and that they are both extremal.

We claim that the two edges labelled $4$ in the above paragraph must 
have an odd number of other edges between them.  
If not, we may reduce to the case
where $n$ is odd and $m(s_1, s_2) = m(s_{n-1}, s_n) = 4$, and now $$
(s_1 s_3 s_5 \cdots s_n)
(s_2 s_4 s_6 \cdots s_{n-1})
(s_1 s_3 s_5 \cdots s_n)
$$ provides the required counterexample.

The parity condition on $n$ now forces
$X = \ti{C}_{n-1}$ for $n$ even, and these graphs are on the list.

We have now reduced to the case where $X$ has at most one edge labelled $4$.
If no such edge exists, we are in type $A$, which is on the list, so suppose
there is a unique edge labelled $4$.  We claim that if this edge is not
an extremal edge (which would give type $B_n$) and not adjacent to an 
extremal edge (which would give type $F_n$), then $X$ must be the
graph of type $\ti{F}_5$ shown in Figure 4.  If not, we may reduce to the
case where $n = 7$ and $m(s_3, s_4) = 4$.  In this case, the required
counterexample can be taken to be $$
(s_3 s_5 s_7)
(s_4 s_6)
(s_3 s_5)
(s_2 s_4)
(s_1 s_3)
(s_2 s_4)
(s_3 s_5)
(s_4 s_6)
(s_3 s_5 s_7)
.$$  Since $\ti{F}_5$ is on the list, our proof is complete.
\qed\enddemo

\head \S\secconc.  Concluding remarks \endhead

Using the techniques of \S\secb, it is possible to derive sharper results
about the structure constants of the $c$-basis for star reducible Coxeter
groups.  In particular, writing $$
c_x c_y = \sum_{w \in W_c} f(x, y, w) c_w
,$$ one may show that all nonzero Laurent polynomials $f(x, y, w)$, for a 
fixed $x$ and $y$, are (positive) integer multiples of the same power of $\d$.

According to \cite{{\bf 4}}, interesting algebras and representations defined over
$\enn$ come from category theory, and are best understood when their 
categorical origin has been discovered.  In \cite{{\bf 9}}, the author showed
how in the case of Coxeter types $A$, $B$, $H$ and $I$, the positivity property
of Theorem \secd.6 (v) may be understood in terms of a category of tangles.
However, there ought to be some representation-theoretic way to understand
this, building on the work of Stroppel \cite{{\bf 20}, \S4} in the case of 
Coxeter types $A$, $B$ and $D$.

\head Acknowledgements \endhead

I thank J. Losonczy for many helpful comments on an early version of this
paper, and for pointing out an error.  I am also grateful to the referee 
for his or her careful reading of the paper and constructive suggestions for
improvements.

\leftheadtext{} \rightheadtext{}
%\vfill\eject

\end